\documentclass[a4paper,oneside,10.5pt]{article}%
\usepackage{makeidx}
\usepackage[english]{babel}
\usepackage{amsmath}
\usepackage{amsfonts}
\usepackage{amssymb}
\usepackage{stmaryrd}
\usepackage{graphicx}
\usepackage{mathrsfs}
\usepackage[colorlinks,linkcolor=red,anchorcolor=blue,citecolor=blue,urlcolor=blue]{hyperref}
\providecommand{\U}[1]{\protect\rule{.1in}{.1in}}
\providecommand{\U}[1]{\protect \rule{.1in}{.1in}}

\newtheorem{theorem}{Theorem}[section]

\newtheorem{lemma}[theorem]{Lemma}

\newtheorem{proposition}[theorem]{Proposition}

\numberwithin{equation}{section}

\begin{document}

\title{On the Seventh Power Moment of $\Delta(x)$}
\author{Jinjiang Li}
\date{}
\maketitle

\begin{center}
  Department of Mathematics, \\
   China University of Mining and Technology, \\
   Beijing 100083, P. R. China \\
   Email: \href{mailto:jinjiang.li.math@gmail.com}{jinjiang.li.math@gmail.com}
\end{center}

{\textbf{Abstract:}} Let $\Delta(x)$ be the error term of the Dirichlet divisor problem. In this paper,
we establish an asymptotic formula of the seventh-power moment of $\Delta(x)$ and prove that
\begin{equation*}
   \int_2^T \Delta^7(x)\mathrm{d}x=
   \frac{7(5s_{7;3}(d)-3s_{7;2}(d)-s_{7;1}(d))}{2816\pi^7}T^{11/4}+O(T^{11/4-\delta_7+\varepsilon})
\end{equation*}
with $\delta_7=1/336,$ which improves the previous result.


\section{Introduction and main result}

    Throughout this paper, let $d(n)$ denote the Dirichlet divisor function. In 1838, Dirichlet proved that
the error term
\begin{equation*}
  \Delta(x):=\sum_{n\leqslant x}d(n)-x\log x-(2\gamma-1)x, \qquad \textrm{with $x\geqslant2,$}
\end{equation*}
satisfies $\Delta(x)\ll x^{1/2}.$ Here $\gamma$ is Euler's constant. Since then the determination of the exact order of $\Delta(x)$ has been called Dirichlet's divisor problem. Many writers have sharpened Dirichlet's bound for $\Delta(x).$  The latest result is due to Huxley \cite{Huxley-2}, who proved that
\begin{equation*}
   \Delta(x)\ll x^{131/416}(\log x)^{26957/8320}.
\end{equation*}
For a survey of the history of this problem, see Kr\"{a}tzel \cite{Kratzel}.

   In the opposite direction, Hardy \cite{Hardy} proved that
\begin{equation*}
   \Delta(x)=\left\{
      \begin{array}{l}
         \Omega_+\left(x^{1/4}(\log x)^{1/4}\log\log x\right), \\
         \Omega_-\left(x^{1/4}\right).
      \end{array}
   \right.
\end{equation*}
The best results in this direction to date are
\begin{equation*}
   \Delta(x)=\Omega_+\left(x^{1/4}(\log x)^{1/4}(\log\log x)^{(3+\log 4)/4}\exp(-c\sqrt{\log\log\log x})\right)
\end{equation*}
and
\begin{equation*}
   \Delta(x)=\Omega_-\left(x^{1/4}\exp(c(\log\log x)^{1/4}(\log\log\log x)^{-3/4}\right)
\end{equation*}
for some constant $c>0,$ due to Hafner \cite{Hafner} and \cite{Corradi} respectively. It is conjectured that
$\Delta(x)\ll x^{1/4+\varepsilon}$ is true for every $\varepsilon>0.$ The evidence in support of this conjecture
has been given by Tong \cite{Tong} and Ivic \cite{Ivic}, who proved, respectively, that
\begin{equation}\label{Tong}
  \int_2^T\Delta^2(x)\mathrm{d}x=\frac{\zeta^4(\frac{3}{2})}{6\pi^2\zeta(3)}T^{3/2}+O(T\log^5T)
\end{equation}
and
\begin{equation}\label{Ivic}
   \int_2^T|\Delta(x)|^A\mathrm{d}x\ll_\varepsilon T^{1+\frac{A}{4}+\varepsilon},\qquad
   \textrm{for $0\leqslant A\leqslant \frac{35}{4}$}
\end{equation}
and any $\varepsilon>0.$ On the other hand , Vorono\"{i} \cite{Voronoi} proved that
\begin{equation}
   \int_{2}^T\Delta(x)\mathrm{d}x=T/4+O(T^{3/4}),
\end{equation}
which in conjunction with (\ref{Tong}) and (\ref{Ivic}) shows that $\Delta(x)$ has a lot of sign changes and cancellations between the positive and negative portions.

   Tsang \cite{Tsang} first studied the third and fourth-power moments of $\Delta(x).$ He proved that
\begin{equation}\label{tsang-3}
   \int_2^T\Delta^3(x)\mathrm{d}x=\frac{3c_1}{28\pi^3}T^{7/4}+O(T^{7/4-\delta_3+\varepsilon}),
\end{equation}
\begin{equation}\label{tsang-4}
   \int_2^T\Delta^4(x)\mathrm{d}x=\frac{3c_2}{64\pi^4}T^{2}+O(T^{2-\delta_4+\varepsilon}),
\end{equation}
where $\delta_3=1/14,\delta_4=1/23$ and
\begin{equation*}
   c_1:= \sum_{\alpha,\beta,h\in\mathbb{N}}(\alpha\beta(\alpha+\beta))^{-3/2}h^{-9/4}|\mu(h)|
             d(\alpha^2h)d(\beta^2h)d((\alpha+\beta)^2h),
\end{equation*}
\begin{equation*}
  c_2:= \sum_{\substack{n,m,k,\ell\in\mathbb{N} \\ \sqrt{n}+\sqrt{m}=\sqrt{k}+\sqrt{\ell}}}
             (nmk\ell)^{-3/4}d(n)d(m)d(k)d(\ell),
\end{equation*}
and $\mu(h)$ is the M\"{o}bius function.

    In \cite{Zhai-1}, Zhai proved that (\ref{tsang-3}) holds for $\delta_3=1/4.$ Ivi\'{c} and Sargos
\cite{Ivic-Sargos} proved that (\ref{tsang-3}) holds for $\delta_3=7/20.$ Following the approach of Tsang \cite{Tsang}, Zhai \cite{Zhai-1} proved that the equation (\ref{tsang-4}) holds for $\delta_4=2/41.$ This approach used the method of exponential sums. In particular, if the exponent pair conjecture is true, namely, if
$(\varepsilon,1/2+\varepsilon)$ is an exponent pair, then the equation (\ref{tsang-4}) holds for $\delta_4=1/14.$
Moreover, in \cite{Ivic-Sargos}, Ivi\'{c} and Sargos proved a substantially better result that the equation (\ref{tsang-4}) holds for $\delta_4=1/12.$ Later, combining the method of \cite{Ivic-Sargos} and a deep result of Robert and Sargos \cite{Robert-Sargos}, Zhai \cite{Zhai-3} proved that the equation (\ref{tsang-4}) holds for $\delta_4=3/28.$ Recently, Kong \cite{Kong} proved that $\delta_4=1/8.$

   By a unified approach, Zhai \cite{Zhai-2} proved that the asymptotic formula
\begin{equation}\label{mean-uniform}
   \int_1^T\Delta^k(x)\mathrm{d}x=C_kT^{1+k/4}+O\left(T^{1+k/4-\delta_k+\varepsilon}\right)
\end{equation}
holds for $3\leqslant k\leqslant9,$ where $C_k$ and $0<\delta_k<1$ are explicit constants. He gives
$\delta_5=1/64, \delta_6=35/4742, \delta_7=17/6312, \delta_8=8/9433, \delta_9=13/75216.$
The asymptotic formula (\ref{mean-uniform}) improved the result of Heath-Brown \cite{Heath-Brown}. When $k=5,$ the
asymptotic formula (\ref{mean-uniform}) holds for $\delta_5=1/64,$ which improved an earlier exponent $\delta_5=5/816$ proved in \cite{Zhai-1} by the approach of Tsang \cite{Tsang}. In \cite{Zhang-deyu}, Zhang and Zhai
improved the previous result of the case $k=5$ and proved $\delta_5=3/80.$ Meanwhile, Wang \cite{Wang} studied the case $k=6$ and proved $\delta_6=3/248,$ which improved the result of Zhai \cite{Zhai-2}, i.e. $\delta_6=35/4742.$

  The aim of this paper is to improve the value of $\delta_7=17/6312,$ which is achieved by Zhai \cite{Zhai-2}. The
main result is the following

\begin{theorem}
    We have
\begin{equation*}
  \int_2^T\Delta^7(x)\mathrm{d}x=
  \frac{7(5s_{7;3}(d)-3s_{7;2}(d)-s_{7;1}(d))}{2816\pi^7}T^{11/4}+O(T^{11/4-\delta_7+\varepsilon})
\end{equation*}
with $\delta_7=1/336,$ where
\begin{eqnarray*}
  & & s_{7;3}(d)=\sum_{\substack{n,m,k,\ell,r,s,q\in\mathbb{N^*}\\
           \sqrt{n}+\sqrt{m}+\sqrt{k}+\sqrt{\ell}=\sqrt{r}+\sqrt{s}+\sqrt{q}}}
          \frac{ d(n)d(m)d(k)d(\ell)d(r)d(s)d(q)} {(nmk\ell rsq)^{3/4}},  \\
  & & s_{7;2}(d)=\sum_{\substack{n,m,k,\ell,r,s,q\in\mathbb{N^*}\\
           \sqrt{n}+\sqrt{m}+\sqrt{k}+\sqrt{\ell}+\sqrt{r}=\sqrt{s}+\sqrt{q}}} \frac{ d(n)d(m)d(k)d(\ell)d(r)d(s)d(q)} {(nmk\ell rsq)^{3/4}},\\
  & & s_{7;1}(d)=\sum_{\substack{n,m,k,\ell,r,s,q\in\mathbb{N^*}\\
           \sqrt{n}+\sqrt{m}+\sqrt{k}+\sqrt{\ell}+\sqrt{r}+\sqrt{s}=\sqrt{q}}}
           \frac{ d(n)d(m)d(k)d(\ell)d(r)d(s)d(q)} {(nmk\ell rsq)^{3/4}}.
\end{eqnarray*}
\end{theorem}

\textbf{Notations.} Throughout this paper, $\|x\|$ denotes the distance from $x$ to the nearest integer, i.e.,
$\|x\|=\min\limits_{n\in\mathbb{Z}}|x-n|.$ $[x]$ denotes the integer part of $x;$ $n\sim N$ means $N<n\leqslant 2N;$ $n\asymp N$ means $C_1N\leqslant n\leqslant C_2N$ with positive constants $C_1,C_2$ satisfying $C_1<C_2.$
$\varepsilon$ always denotes an arbitrary  small positive constant which may not be the same at different occurances. We shall use the estimates $d(n)\ll n^{\varepsilon}.$ Suppose $f: \mathbb{N}\to \mathbb{R}$ is any function satisfying $f(n)\ll n^{\varepsilon},$ $k\geqslant2$ is a fixed integer. Define
\begin{equation}\label{series}
   s_{k;\ell}(f):=\sum_{\substack{n_1,\cdots,n_\ell,n_{\ell+1},\cdots,n_k\in\mathbb{N^*}\\
   \sqrt{n_1}+\cdots+\sqrt{n_{\ell}}=\sqrt{n_{\ell+1}}+\cdots+\sqrt{n_k} }}
   \frac{f(n_1)f(n_2)\cdots f(n_k)}{(n_1n_2\cdots n_k)^{3/4}}, \quad 1\leqslant\ell <k.
\end{equation}
We shall use $s_{k;\ell}(f)$ to denote both of the series (\ref{series}) and its value. Suppose $y>1$ is a large
parameter, and we define
\begin{equation*}
   s_{k;\ell}(f;y):=\sum_{\substack{n_1,\cdots,n_\ell,n_{\ell+1},\cdots,n_k\leqslant y\\
   \sqrt{n_1}+\cdots+\sqrt{n_{\ell}}=\sqrt{n_{\ell+1}}+\cdots+\sqrt{n_k} }}
   \frac{f(n_1)f(n_2)\cdots f(n_k)}{(n_1n_2\cdots n_k)^{3/4}}, \quad 1\leqslant\ell <k.
\end{equation*}

\section{Preliminary Lemmas}

\begin{lemma}\label{lemma-1}
      Suppose $k\geqslant3, (i_1,i_2,\cdots,i_{k-1})\in\{0,1\}^{k-1}$ such that
   \begin{equation*}
     \sqrt{n_1}+(-1)^{i_1}\sqrt{n_2}+(-1)^{i_2}\sqrt{n_3}+\cdots+(-1)^{i_{k-1}}\sqrt{n_k}\neq0.
   \end{equation*}
 Then we have
 \begin{equation*}
     |\sqrt{n_1}+(-1)^{i_1}\sqrt{n_2}+(-1)^{i_2}\sqrt{n_3}+\cdots+(-1)^{i_{k-1}}\sqrt{n_k}|
     \gg\max(n_1,n_2,\cdots,n_k)^{-(2^{k-2}-2^{-1})}.
   \end{equation*}
\end{lemma}
\textbf{Proof.} See Lemma $2.2$ of \cite{Zhai-2}.

\begin{lemma} \label{lemma-2}
If $g(x)$ and $h(x)$ are continuous real-valued functions of $x$ and $g(x)$ is monotonic, then
  \begin{equation*}
   \int_a^bg(x)h(x)\mathrm{d}x\ll\left(\max_{a\leqslant x\leqslant b}|g(x)|\right)
   \left(\max_{a\leqslant u<v\leqslant b}\left|\int_u^vh(x)\mathrm{d}x\right|\right).
  \end{equation*}
\end{lemma}
\textbf{Proof.} See Lemma $1$ of \cite{Tsang}.

\begin{lemma} \label{lemma-3}
Suppose $A,B\in\mathbb{R}, A\neq0.$ Then
  \begin{equation*}
      \int_T^{2T}\cos(A\sqrt{t}+B)\mathrm{d}t\ll T^{1/2}|A|^{-1}.
  \end{equation*}
\end{lemma}
\textbf{Proof.} It follows from Lemma \ref{lemma-2} easily.

\begin{lemma}\label{lemma-4}
   Suppose $K\geqslant10,\alpha,\beta\in\mathbb{R},2K^{-1/2}\leqslant|\alpha|\ll K^{1/2}$ and $0<\delta<1/2.$
   Then we have
   \begin{equation*}
     \#\{k\sim K: \|\beta+\alpha\sqrt{t}\|<\delta\}\ll K\delta+K^{1/2+\varepsilon}.
   \end{equation*}
\end{lemma}
\textbf{Proof.} See Lemma $4$ of \cite{Zhai-3}.

\begin{lemma}\label{lemma-5}
   Let $a,\delta$ be real numbers, $0<\delta<a/4,$ and let $k$ be a positive integer. There exists a function $\varphi(y)$
   which is $k$ times continuously differentiable and such that

  \begin{equation*}
    \left\{
      \begin{array}{cll}
          \varphi(y)=1,    & &  \textrm{for \quad} |y|\leqslant a-\delta, \\
          0<\varphi(y)<1,  & &  \textrm{for \quad} a-\delta<|y|< a+\delta, \\
          \varphi(y)=0,    & &  \textrm{for \quad} |y|\geqslant a+\delta,
      \end{array}
    \right.
  \end{equation*}
  and its Fourier transform
   \begin{equation*}
      \Phi(x)=\int_{-\infty}^{+\infty} e(-xy)\varphi(y)\mathrm{d}y
   \end{equation*}
   satisfies the inequality
   \begin{equation*}
      \left|\Phi(x)\right|\leqslant\min\left(2a,\frac{1}{\pi|x|},\frac{1}{\pi|x|}\left(\frac{k}{2\pi|x|\delta}\right)^k\right).
   \end{equation*}
\end{lemma}
\textbf{Proof.} See  \cite{Piatetski-Shapiro} or \cite{Segal}.

\begin{lemma}\label{lemma-6}
    Let $d(n)$ denote the divisor function. Then we have
   \begin{equation*}
      |s_{k;\ell}(d)-s_{k;\ell}(d;y)|\ll y^{-1/2+\varepsilon},\qquad 1\leqslant\ell<k.
   \end{equation*}
\end{lemma}
\textbf{Proof.} See Lemma $3.1$ of \cite{Zhai-2}.

\begin{lemma} \label{lemma-7}
   Suppose $k\geqslant3; (i_1,i_2,\cdots,i_{k-1})\in\{0,1\}^{k-1}, (i_1,i_2,\cdots,i_{k-1})\neq(0,0,\cdots,0);\\
1<N_1,N_2,\cdots,N_k; 0<\Delta\ll E^{1/2}, E=\max(N_1,N_2,\cdots,N_k).$ Let
   \begin{equation*}
      \mathscr{A}=\mathscr{A}(N_1,N_2,\cdots,N_k;i_1,i_2,\cdots,i_k;\Delta)
   \end{equation*}
   denote the number of solutions of the inequality
   \begin{equation*}
     |\sqrt{n_1}+(-1)^{i_1}\sqrt{n_2}+(-1)^{i_2}\sqrt{n_3}+\cdots+(-1)^{i_{k-1}}\sqrt{n_k}|<\Delta
   \end{equation*}
   with $n_j\sim N_j, j=1,2,\cdots,k.$ Then
   \begin{equation*}
       \mathscr{A}\ll \Delta E^{-1/2}N_1N_2\cdots N_k+E^{-1}N_1N_2\cdots N_k.
   \end{equation*}
\end{lemma}
\textbf{Proof.} See Lemma $2.4$ of \cite{Zhai-2}.

\begin{lemma} \label{lemma-8}
   Suppose $1\leqslant N\leqslant M\leqslant K\leqslant L, 1\leqslant R\leqslant S\leqslant Q, L\asymp Q$ and
   $0<\Delta\ll Q^{1/2}.$ Let $\mathscr{A}_1(N,M,K,L,R,S,Q;\Delta)$ denote the number of solutions of the inequality
   \begin{equation}\label{lemma-8-1}
      0<|\sqrt{n}+\sqrt{m}+\sqrt{k}+\sqrt{\ell}-\sqrt{r}-\sqrt{s}-\sqrt{q}|<\Delta
   \end{equation}
   with $n\sim N, m\sim M, k\sim K, \ell\sim L, r\sim R, s\sim S, q\sim Q.$ Then
   \begin{equation*}
     \mathscr{A}_1(N,M,K,L,R,S,Q;\Delta)\ll \Delta Q^{1/2}NMKLRS+NMKRSL^{1/2+\varepsilon}.
   \end{equation*}
   In particular, if $\Delta Q^{1/2}\gg1,$ then
   \begin{equation*}
     \mathscr{A}_1(N,M,K,L,R,S,Q;\Delta)\ll \Delta Q^{1/2}NMKLRS.
   \end{equation*}
\end{lemma}

 \textbf{Proof.} If $(n,m,k,\ell.r,s,q)$ satisfies (\ref{lemma-8-1}), then
  \begin{equation*}
     (\sqrt{n}+\sqrt{m}+\sqrt{k}-\sqrt{r}-\sqrt{s})+\sqrt{\ell}=\sqrt{q}+\theta\Delta
  \end{equation*}
for some $0<|\theta|<1.$ Thus, we have
\begin{equation*}
  2\ell^{1/2}(\sqrt{n}+\sqrt{m}+\sqrt{k}-\sqrt{r}-\sqrt{s})+(\sqrt{n}+\sqrt{m}+\sqrt{k}-\sqrt{r}-\sqrt{s})^2+\ell
  =q+u
\end{equation*}
with $|u|=|2q^{1/2}\theta\Delta+\theta^2\Delta^2|\leqslant2q^{1/2}\Delta+\Delta^2\ll\Delta Q^{1/2}.$ Then we have
\begin{equation*}
  q=2\ell^{1/2}(\sqrt{n}+\sqrt{m}+\sqrt{k}-\sqrt{r}-\sqrt{s})+(\sqrt{n}+\sqrt{m}+\sqrt{k}-\sqrt{r}-\sqrt{s})^2+\ell-u
\end{equation*}
with $|u|\leqslant C\Delta Q^{1/2}$ for some absolute positive constant $C>0.$ Hence the quantity of \\
$\mathscr{A}_1(N,M,K,L,R,S,Q;\Delta)$ does not exceed the number of solutions of
\begin{equation}\label{Lemma-8-2}
 \left|2\ell^{1/2}(\sqrt{n}+\sqrt{m}+\sqrt{k}-\sqrt{r}-\sqrt{s})+
 (\sqrt{n}+\sqrt{m}+\sqrt{k}-\sqrt{r}-\sqrt{s})^2+\ell-q\right|
 <C\Delta Q^{1/2}
\end{equation}
 with $n\sim N, m\sim M, k\sim K, \ell\sim L, r\sim R, s\sim S, q\sim Q.$

  If $\Delta Q^{1/2}\gg1,$ then for fixed $(n,m,k,\ell,r,s),$ the number of $q$ for which (\ref{Lemma-8-2}) holds
is $\ll1+\Delta Q^{1/2}\ll\Delta Q^{1/2}.$ Hence
\begin{equation*}
     \mathscr{A}_1(N,M,K,L,R,S,Q;\Delta)\ll \Delta Q^{1/2}NMKLRS.
   \end{equation*}

  Now suppose $\Delta Q^{1/2}\leqslant1/4C.$ Then for fixed $(n,m,k,\ell,r,s),$ there is at most one $q$ such that
(\ref{Lemma-8-2}) holds. If such $q$ exists, then we have
  \begin{equation}\label{lemma-8-3}
     \left\|2\ell^{1/2}(\sqrt{n}+\sqrt{m}+\sqrt{k}-\sqrt{r}-\sqrt{s})+
     (\sqrt{n}+\sqrt{m}+\sqrt{k}-\sqrt{r}-\sqrt{s})^2\right\|<C\Delta Q^{1/2}.
  \end{equation}
We shall use Lemma \ref{lemma-4} to bound the number of solutions of (\ref{lemma-8-3}) with
$\alpha=2(\sqrt{n}+\sqrt{m}+\sqrt{k}-\sqrt{r}-\sqrt{s}),\beta=(\sqrt{n}+\sqrt{m}+\sqrt{k}-\sqrt{r}-\sqrt{s})^2.$
Let $\mathscr{D}_1$ denote the number of solutions of (\ref{lemma-8-3}) with $|\alpha|\geqslant2L^{-1/2},$ and
$\mathscr{D}_2$ the number of solutions of (\ref{lemma-8-3}) with $|\alpha|<2L^{-1/2}.$ By Lemma \ref{lemma-4},
we get
\begin{eqnarray*}
  \mathscr{D}_1 & \ll & (\Delta Q^{1/2}L+L^{1/2+\varepsilon})NMKRS \\
                & \ll & \Delta Q^{1/2}NMKLRS+NMKRSL^{1/2+\varepsilon}.
\end{eqnarray*}
Now we estimate $\mathscr{D}_2.$ From $|\alpha|<2L^{-1/2},$ we can get $K\asymp R.$ If
$\sqrt{n}+\sqrt{m}+\sqrt{k}=\sqrt{r}+\sqrt{s},$ then from (\ref{Lemma-8-2}), we get $\ell=q.$ This contradicts to the fact that $|\sqrt{n}+\sqrt{m}+\sqrt{k}+\sqrt{\ell}-\sqrt{r}-\sqrt{s}\sqrt{q}|>0.$ Therefore, we have
$\sqrt{n}+\sqrt{m}+\sqrt{k}\neq\sqrt{r}+\sqrt{s}.$ By Lemma \ref{lemma-1}, we have
$|\sqrt{n}+\sqrt{m}+\sqrt{k}-\sqrt{r}+\sqrt{s}|\gg S^{-15/2}$ for any such $(n,m,k,r,s).$ By a splitting argument and Lemma \ref{lemma-7}, there exists a $\delta$ satisfying $S^{-15/2}\ll\delta\ll L^{-1/2},$ which holds
\begin{eqnarray*}
   \mathscr{D}_2 & \ll & \log2S\cdot \sum_{\delta<|\sqrt{n}+\sqrt{m}+\sqrt{k}-\sqrt{r}+\sqrt{s}|\leqslant2\delta}1\\
   & \ll & \log2S\cdot (\delta S^{1/2}NMKR+NMKR) \\
   & \ll & L^{-1/2}S^{1/2+\varepsilon}NMKR+NMKRS^{\varepsilon} \\
   & \ll & NMKRS^{\varepsilon},
\end{eqnarray*}
which can be absorbed into the estimate of $\mathscr{D}_1.$ This completes the proof of Lemma \ref{lemma-8}.

\begin{lemma} \label{lemma-9}
     Suppose $1\leqslant N\leqslant M\leqslant K\leqslant L, 1\leqslant R\leqslant S\leqslant Q, L\asymp Q$ and
   $0<\Delta\ll Q^{1/2}.$ Let $\mathscr{A}_{2,\pm}(N,M,K,L,R,S,Q;\Delta)$ denote the number of solutions of the inequality
   \begin{equation*}
      0<|\sqrt{n}+\sqrt{m}+\sqrt{k}+\sqrt{\ell}+\sqrt{r}\pm\sqrt{s}-\sqrt{q}|<\Delta
   \end{equation*}
   with $n\sim N, m\sim M, k\sim K, \ell\sim L, r\sim R, s\sim S, q\sim Q.$ Then
   \begin{equation*}
     \mathscr{A}_{2,\pm}(N,M,K,L,R,S,Q;\Delta)\ll \Delta Q^{1/2}NMKLRS+NMKRSL^{1/2+\varepsilon}.
   \end{equation*}
   In particular, if $\Delta Q^{1/2}\gg1,$ then
   \begin{equation*}
     \mathscr{A}_{2,\pm}(N,M,K,L,R,S,Q;\Delta)\ll \Delta Q^{1/2}NMKLRS.
   \end{equation*}
\end{lemma}
\textbf{Proof.} The proof of Lemma \ref{lemma-9} is similar to that of Lemma \ref{lemma-8}, so we omit the details.

\begin{lemma} \label{lemma-10}
   Suppose $N_j\geqslant2\quad(j=1,2,3,4,5,6,7)$ are real numbers, $\Delta>0.$ let \\
$\mathscr{A}_{\pm,\pm}(N_1.N_2,N_3,N_4,N_5,N_6,N_7;\Delta)$ denote the number of solutions of the inequality
\begin{equation*}
  0<|\sqrt{n_1}+\sqrt{n_2}+\sqrt{n_3}+\sqrt{n_4}\pm\sqrt{n_5}\pm\sqrt{n_6}-\sqrt{n_7}|<\Delta
\end{equation*}
with $n_j\sim N_j,(j=1,2,3,4,5,6,7),n_j\in\mathbb{N}^*.$ Then we have
\begin{equation*}
  \mathscr{A}_{\pm,\pm}(N_1.N_2,N_3,N_4,N_5,N_6,N_7;\Delta)\ll\prod_{j=1}^7
\left(\Delta^{1/7}N_j^{13/14}+N_j^{5/7}\right)N_j^{\varepsilon}.
\end{equation*}
\end{lemma}
\textbf{Proof.} Taking $a=6\Delta/5,\delta=\Delta/5$ in Lemma \ref{lemma-5}, there exists a function $\varphi_1(y),$
 which is $\ell=[7\log(N_1N_2N_3N_4N_5N_6N_7)]$ times continuously differentiable such that
\begin{equation*}
    \left\{
      \begin{array}{cll}
          \varphi_1(y)=1,    & &  \textrm{if \quad} |y|\leqslant \Delta, \\
          0<\varphi_1(y)<1,  & &  \textrm{for \quad} \Delta<|y|< 7\Delta/5, \\
          \varphi_1(y)=0,    & &  \textrm{for \quad} |y|\geqslant 7\Delta/5.
      \end{array}
    \right.
  \end{equation*}
Let
   \begin{equation*}
      \Phi_1(x)=\int_{-\infty}^{+\infty} e(-xy)\varphi_1(y)\mathrm{d}y,
   \end{equation*}
 then it satisfies
   \begin{equation}\label{lemma-10-1}
      \left|\Phi_1(x)\right|\leqslant
\min\left(\frac{12\Delta}{5},\frac{1}{\pi|x|},\frac{1}{\pi|x|}\left(\frac{5\ell}{2\pi|x|\Delta}\right)^\ell\right),
   \end{equation}
and
\begin{equation}\label{lemma-10-2}
  \varphi_1(y)=\int_{-\infty}^{+\infty}e(xy)\Phi_1(x)\mathrm{d}x.
\end{equation}
Set
 \begin{equation*}
   R_{\pm,\pm}=\sum_{\substack{n_j\sim N_j \\j=1,2,3,4,5,6,7}}
\varphi_1(\sqrt{n_1}+\sqrt{n_2}+\sqrt{n_3}+\sqrt{n_4}\pm\sqrt{n_5}\pm\sqrt{n_6}-\sqrt{n_7}).
 \end{equation*}
By the definition of $\varphi_1(y),$ we get
 \begin{equation}\label{lemma-10-3}
    \mathscr{A}_{\pm,\pm}(N_1.N_2,N_3,N_4,N_5,N_6,N_7;\Delta)\leqslant R_{\pm,\pm}.
 \end{equation}
 We estimate $R_{-,-}$ first. By (\ref{lemma-10-2}), we have
\begin{eqnarray*}
  R_{-,-} & =& \sum_{\substack{n_j\sim N_j \\j=1,2,3,4,5,6,7}}
               \varphi_1(\sqrt{n_1}+\sqrt{n_2}+\sqrt{n_3}+\sqrt{n_4}-\sqrt{n_5}-\sqrt{n_6}-\sqrt{n_7}) \\
   & = & \sum_{\substack{n_j\sim N_j \\j=1,2,3,4,5,6,7}}
        \int_{-\infty}^{+\infty}
        e\left(x(\sqrt{n_1}+\sqrt{n_2}+\sqrt{n_3}+\sqrt{n_4}-\sqrt{n_5}-\sqrt{n_6}-\sqrt{n_7})\right)
         \Phi_1(x)\mathrm{d}x.
\end{eqnarray*}
Let $S(x;N):=\sum\limits_{n\sim N}e(x\sqrt{n}),$ we have
\begin{equation*}
  R_{-,-}=\int_{-\infty}^{+\infty}
          S(x;N_1)S(x;N_2)S(x;N_3)S(x;N_4)\overline{S(x;N_5)S(x;N_6)S(x;N_7)} \Phi_1(x)\mathrm{d}x,
\end{equation*}
if we notice that $\overline{S(x;N)}=\sum\limits_{n\sim N}e(-x\sqrt{n}).$ Applying H\"{o}lder¡¯s inequality, we get
\begin{equation}\label{lemma-10-4}
   R_{-,-}\leqslant\prod_{j=1}^7\left(\int_{-\infty}^{+\infty} |S(x;N_j)|^7|\Phi_1(x)|\mathrm{d}x\right)^{1/7}.
\end{equation}
Let
\begin{equation*}
  T(N):=\int_{0}^{+\infty} |S(x;N)|^7|\Phi_1(x)|\mathrm{d}x.
\end{equation*}
It is sufficient to estimate $T(N),$ where $N=N_j$ for some $j\in\{1,2,3,4,5,6,7\}.$ Let
\begin{equation*}
   K:=\frac{100\ell^2}{\Delta},\quad\ell=[7\log(N_1N_2N_3N_4N_5N_6N_7)];\quad K_0:=N^{1/2}.
\end{equation*}
Using the trivial estimate $S(x;N)\ll N$ and the estimate
   \begin{equation*}
      \left|\Phi_1(x)\right|\leqslant\frac{1}{\pi|x|}\left(\frac{5\ell}{2\pi|x|\Delta}\right)^\ell,
   \end{equation*}
we have
\begin{eqnarray}\label{lemma-10-5}
          \int_{K}^\infty|S(x;N)|^7|\Phi_1(x)|\mathrm{d}x
 & \ll &  N^7\int_{K}^\infty|\Phi_1(x)|\mathrm{d}x   \ll
          N^7\left(\frac{5\ell}{2\pi\Delta}\right)^\ell\int_{K}^\infty\frac{1}{x^{\ell+1}}\mathrm{d}x
                 \nonumber  \\
 & \ll &  \frac{N^75^\ell}{\ell} \left(\frac{\ell}{2\pi K\Delta}\right)^\ell \ll
           \frac{N^75^\ell}{\ell}\frac{1}{\ell^\ell}
                 \nonumber  \\
 & \ll & \frac{N^7(N_1\cdots N_7)^{7\log5}}
            {(N_1\cdots N_7)^{7\log7}(N_1\cdots N_7)^{7\log\log(N_1\cdots N_7)}} \ll1.
\end{eqnarray}
For the mean square of $S(x;N),$ we have
\begin{eqnarray*}
   \int_0^{K_0} |S(x;N)|^2\mathrm{d}x
       & = & \int_0^{K_0}\sum_{n\sim N}\sum_{m\sim N} e(x(\sqrt{n}-\sqrt{m}))\mathrm{d}x \\
       & = & \sum_{n\sim N}\sum_{m\sim N}\int_0^{K_0} e(x(\sqrt{n}-\sqrt{m}))\mathrm{d}x  \\
    & = & \bigg\{\mathop{\sum_{n\sim N}\sum_{m\sim N}}_{n=m}+\mathop{\sum_{n\sim N}\sum_{m\sim N}}_{n\neq m}\bigg\}
          \int_0^{K_0} e(x(\sqrt{n}-\sqrt{m}))\mathrm{d}x.
\end{eqnarray*}
If we use the trivial estimate $S(x;N)\ll N,$ then
\begin{equation*}
    \mathop{\sum_{n\sim N}\sum_{m\sim N}}_{n=m} \int_0^{K_0} e(x(\sqrt{n}-\sqrt{m}))\mathrm{d}x\ll N^{3/2}.
\end{equation*}
For the case $n\neq m,$ we have
\begin{eqnarray*}
 &    & \mathop{\sum_{n\sim N}\sum_{m\sim N}}_{n\neq m} \int_0^{K_0} e(x(\sqrt{n}-\sqrt{m}))\mathrm{d}x\\
 & \ll & \mathop{\sum_{n\sim N}\sum_{m\sim N}}_{n\neq m} \frac{1}{|\sqrt{n}-\sqrt{m}|}
  \ll N^{1/2}\mathop{\sum_{n\sim N}\sum_{m\sim N}}_{n\neq m}  \frac{1}{|n-m|} \\
 & \ll & N^{1/2} \sum_{n\sim N}\sum_{1\leqslant r\ll N}\frac{1}{r}
   \ll   N^{3/2}\log N.
\end{eqnarray*}
Therefore, we have
\begin{equation}\label{lemma-10-6}
 \int_0^{K_0}|S(x;N)|^2\mathrm{d}x\ll N^{3/2}\log N.
\end{equation}
If we notice $|\Phi_1(x)|\ll \Delta$ from (\ref{lemma-10-1}) and the trivial estimate $S(x;N)\ll N,$ then
\begin{eqnarray}\label{lemma-10-7}
   &     &  \int_0^{K_0}|S(x;N)|^7|\Phi_1(x)|\mathrm{d}x \nonumber \\
   & \ll &  \Delta N^5\int_0^{K_0}|S(x;N)|^2\mathrm{d}x \nonumber \\
   & \ll &  \Delta N^{13/2}\log N.
\end{eqnarray}
If $K\leqslant K_0,$ then from (\ref{lemma-10-5}) and (\ref{lemma-10-7}) we get
\begin{equation}\label{lemma-10-8}
  T(N)\ll \Delta N^{13/2}\log N.
\end{equation}
Now suppose $K_0<K.$ By a splitting argument, we have
\begin{equation}\label{lemma-10-9}
 \int_{K_0}^K |S(x;N)|^7|\Phi_1(x)|\mathrm{d}x\ll\Delta\log K\times
\max_{K_0\leqslant U\leqslant K}\int_U^{2U}|S(x;N)|^7\mathrm{d}x.
\end{equation}
From Lemma \ref{lemma-1}, we can get $\Delta^{-1}\ll\max(N_1,N_2,N_3,N_4,N_5,N_6,N_7)^{63/2}$ and thus
$\log K\ll\ell^2.$ On the other hand, we have

\begin{eqnarray}\label{lemma-10-10}
            \int_U^{2U}|S(x;N)|^7\mathrm{d}x
    & \ll & \max_{U\leqslant x\leqslant2U}|S(x;N)|^2\times\int_U^{2U}|S(x;N)|^5\mathrm{d}x
                \nonumber \\
    & \ll &  N^2(N^{9/2}+UN^3)N^\varepsilon
                \nonumber \\
    & \ll & (N^{13/2}+UN^5)N^\varepsilon ,
\end{eqnarray}
which is derived from equation (2.18) of Zhang and Zhai \cite{Zhang-deyu}.

 From (\ref{lemma-10-9}) and (\ref{lemma-10-10}) and noticing $\Delta K\ll\ell^2,$ we get
\begin{eqnarray*}
 \int_{K_0}^K |S(x;N)|^7|\Phi_1(x)|\mathrm{d}x
      & \ll &  \Delta\log K\times (N^{13/2}+UN^5)N^\varepsilon   \\
      & \ll &  (\Delta N^{13/2}+N^5)N^\varepsilon,
\end{eqnarray*}
which combining (\ref{lemma-10-5}) and (\ref{lemma-10-7}) gives
\begin{equation}\label{lemma-10-11}
   T(N)=\int_0^\infty|S(x;N)|^7|\Phi_1(x)|\mathrm{d}x\ll(\Delta N^{13/2}+N^5)N^\varepsilon.
\end{equation}

From (\ref{lemma-10-3}), (\ref{lemma-10-4}) and (\ref{lemma-10-11}), we get the result of Lemma \ref{lemma-10}
for the case ``$-,-$''. By noting the properties of conjugation, the estimates of other cases are exactly the
same as that of the case ``$-,-$''. This completes the proof of Lemma \ref{lemma-10}.

\section{Proof of Theorem}

  In this section, we shall prove the theorem. We begin with the following truncated form of the Vorono\"{i}'s formula (\cite{Ivic-book}, equation (2.25) ), i.e.
\begin{equation}\label{voronoi}
 \Delta(x)=\frac{1}{\sqrt{2}\pi}\sum_{n\leqslant N}\frac{d(n)}{n^{3/4}}x^{1/4}\cos(4\pi\sqrt{nx}-\pi/4)
+O(x^{1/2+\varepsilon}N^{-1/2}),
\end{equation}
for  $1\leqslant N\ll x.$ Set $\Delta(x):=R_1+R_2,$
where
\begin{equation*}
  R_1:=R_1(x)=\frac{1}{\sqrt{2}\pi}\sum_{n\leqslant y}\frac{d(n)}{n^{3/4}}x^{1/4}\cos(4\pi\sqrt{nx}-\pi/4),
  \quad R_2:=R_2(x)=\Delta(x)-R_1.
\end{equation*}
Take $y=T^{1/4}.$ By the elementary estimate $(a+b)^7-a^7\ll |b|a^6+|b|^7,$ we have

\begin{equation*}
 \int_1^T\Delta^7(x)\mathrm{d}x=\int_1^T R_1^7(x)\mathrm{d}x+
O\left(\int_{1}^T|R_1|^6|R_2|\mathrm{d}x+\int_{1}^T|R_2|^7\mathrm{d}x\right).
\end{equation*}
By a splitting argument, it is sufficient to prove the result in the interval $[T,2T].$ We will divide the
process of the proof of the theorem into two parts.

 \begin{proposition}\label{mingti-1}
    For fixed $T\geqslant10, N=T, y=T^{1/4},$ we have
   \begin{equation}\label{proposition-1}
    \int_T^{2T}R_1^7\mathrm{d}x=\frac{7(5s_{7;3}(d)-3s_{7;2}(d)-s_{7;1}(d))}{2816\pi^7}T^{11/4}+
        O(T^{11/4-1/336+\varepsilon}).
   \end{equation}
 \end{proposition}
\textbf{Proof.} Let
\begin{equation*}
 g:=g(n,m,k,\ell,r,s,q)=\left\{
     \begin{array}{cl}
       \frac{d(n)d(m)d(k)d(\ell)d(r)d(s)d(q)}{(nmk\ell rsq)^{3/4}},
                             & \textrm{if\quad} n,m,k,\ell,r,s,q\leqslant y,\\
        0, & \textrm{otherwise.}
     \end{array}
   \right.
\end{equation*}
According to the elementary formula
\begin{equation*}
  \cos a_1\cos a_2\cdots\cos a_h=\frac{1}{2^{h-1}}\sum_{(i_1,i_2,\cdots,i_{h-1})\in\{0,1\}^{h-1}}
   \cos\left(a_1+(-1)^{i_1}a_2+\cdots+(-1)^{i_{h-1}}a_h\right),
\end{equation*}
we can write
\begin{equation*}
  R_1^7=S_1(x)+S_2(x)+S_3(x)+S_4(x)+S_5(x)+S_6(x)+S_7(x),
\end{equation*}
where
\begin{eqnarray*}
 S_1(x) & := & \frac{35}{64}\cos\frac{\pi}{4}
               \sum_{\substack{n,m,k,\ell,r,s,q\leqslant y \\ \sqrt{n}+\sqrt{m}+\sqrt{k}+\sqrt{\ell}
                =\sqrt{r}+\sqrt{s}+\sqrt{q}}} gx^{7/4},
                     \\
  S_2(x) & := & \frac{35}{64} \sum_{\substack{n,m,k,\ell,r,s,q\leqslant y \\ \sqrt{n}+\sqrt{m}+\sqrt{k}+\sqrt{\ell}
                \neq\sqrt{r}+\sqrt{s}+\sqrt{q}}} gx^{7/4}
                     \\
        &  & \times     \cos\left(4\pi\left(\sqrt{n}+\sqrt{m}+\sqrt{k}+\sqrt{\ell}-\sqrt{r}-\sqrt{s}-\sqrt{q}\right)
                     \sqrt{x}-\frac{\pi}{4}\right),
                          \\
   S_3(x) & := & \frac{21}{64}\cos\frac{3\pi}{4}
                 \sum_{\substack{n,m,k,\ell,r,s,q\leqslant y \\ \sqrt{n}+\sqrt{m}+\sqrt{k}+\sqrt{\ell}
                +\sqrt{r}=\sqrt{s}+\sqrt{q}}} gx^{7/4},
                           \\
S_4(x) & := & \frac{21}{64}
                 \sum_{\substack{n,m,k,\ell,r,s,q\leqslant y \\ \sqrt{n}+\sqrt{m}+\sqrt{k}+\sqrt{\ell}
                     +\sqrt{r}\neq\sqrt{s}+\sqrt{q}}}  gx^{7/4}
                             \\
      &  & \times
                 \cos\left(4\pi\left(\sqrt{n}+\sqrt{m}+\sqrt{k}+\sqrt{\ell}+\sqrt{r}-\sqrt{s}-\sqrt{q}\right)
                     \sqrt{x}-\frac{3\pi}{4}\right),
\end{eqnarray*}
\begin{eqnarray*}
 S_5(x) & := & \frac{7}{64} \cos\frac{5\pi}{4}
                 \sum_{\substack{n,m,k,\ell,r,s,q\leqslant y \\ \sqrt{n}+\sqrt{m}+\sqrt{k}+\sqrt{\ell}
                 +\sqrt{r}+\sqrt{s}=\sqrt{q}}} gx^{7/4},
                             \\
   S_6(x) & := & \frac{7}{64}
                 \sum_{\substack{n,m,k,\ell,r,s,q\leqslant y \\ \sqrt{n}+\sqrt{m}+\sqrt{k}+\sqrt{\ell}
                 +\sqrt{r}+\sqrt{s}\neq\sqrt{q}}} gx^{7/4}
                             \\
     &  & \times \cos\left(4\pi\left(\sqrt{n}+\sqrt{m}+\sqrt{k}+\sqrt{\ell}+\sqrt{r}+\sqrt{s}-\sqrt{q}\right)
                     \sqrt{x}-\frac{5\pi}{4}\right),
                             \\
   S_7(x) & := & \frac{1}{64}
                 \sum_{\substack{n,m,k,\ell,r,s,q\leqslant y }} gx^{7/4}
                             \\
    &  & \times \cos\left(4\pi\left(\sqrt{n}+\sqrt{m}+\sqrt{k}+\sqrt{\ell}+\sqrt{r}+\sqrt{s}+\sqrt{q}\right)
                     \sqrt{x}-\frac{7\pi}{4}\right).
\end{eqnarray*}
By Lemma \ref{lemma-6}, we get
\begin{eqnarray} \label{proposition-2}
  \int_T^{2T}S_1(x)\mathrm{d}x & = & \frac{35\sqrt{2}}{128}s_{7;3}(d;y)\int_T^{2T}x^{7/4}\mathrm{d}x
                       \nonumber \\
   & = & \frac{35\sqrt{2}}{128}s_{7;3}(d)\int_T^{2T}x^{7/4}\mathrm{d}x+O\left(T^{11/4}y^{-1/2+\varepsilon}\right)
                        \nonumber \\
   & = & \frac{35\sqrt{2}}{128}s_{7;3}(d)\int_T^{2T}x^{7/4}\mathrm{d}x+O\left(T^{11/4-1/8+\varepsilon}\right).
\end{eqnarray}
Similarly, we can get
\begin{equation}\label{proposition-3}
   \int_T^{2T}S_3(x)\mathrm{d}x=-\frac{21\sqrt{2}}{128}s_{7;2}(d)
        \int_T^{2T}x^{7/4}\mathrm{d}x+O\left(T^{11/4-1/8+\varepsilon}\right)
\end{equation}
and
\begin{equation}\label{proposition-4}
   \int_T^{2T}S_5(x)\mathrm{d}x=-\frac{7\sqrt{2}}{128}s_{7;1}(d)
           \int_T^{2T}x^{7/4}\mathrm{d}x+O\left(T^{11/4-1/8+\varepsilon}\right).
\end{equation}
We now proceed to consider the contribution of $S_7(x).$ Applying Lemma \ref{lemma-3}, we have
\begin{eqnarray}\label{proposition-5}
  \int_T^{2T}S_7(x)\mathrm{d}x & \ll & \sum_{n,m,k,\ell,r,s,q\leqslant y}
    \frac{gT^{9/4}}{\sqrt{n}+\sqrt{m}+\sqrt{k}+\sqrt{\ell}+\sqrt{r}+\sqrt{s}+\sqrt{q}}
                \nonumber \\
  & \ll & T^{9/4+\varepsilon}\sum_{n\leqslant m\leqslant k\leqslant\ell \leqslant r
           \leqslant s\leqslant q\leqslant y}\frac{1}{(nmk\ell rsq)^{3/4}q^{1/2}}
                  \nonumber \\
  & \ll & T^{9/4+\varepsilon}y^{5/4}
   \ll    T^{41/16+\varepsilon}.
\end{eqnarray}
Now we consider the contribution of $S_2(x).$ By the first derivative test, we get
\begin{eqnarray}\label{proposition-6}
  \int_T^{2T}S_2(x)\mathrm{d}x & \ll &
       \sum_{\substack{n,m,k,\ell,r,s,q\leqslant y \\ \sqrt{n}+\sqrt{m}+\sqrt{k}+\sqrt{\ell}
                \neq\sqrt{r}+\sqrt{s}+\sqrt{q}}} g
                       \nonumber \\
    & & \times  \min\left(T^{11/4},\frac{T^{9/4}}{|\sqrt{n}+\sqrt{m}+\sqrt{k}+\sqrt{\ell}
                       -\sqrt{r}-\sqrt{s}-\sqrt{q}| }\right)
                           \nonumber \\
   & \ll & x^{\varepsilon}G(N,M,K,L,R,S,Q),
\end{eqnarray}
where
\begin{eqnarray*}
  G(N,M,K,L,R,S,Q) & = &\sum_{\substack{ \sqrt{n}+\sqrt{m}+\sqrt{k}+\sqrt{\ell}\neq\sqrt{r}+\sqrt{s}+\sqrt{q}\\
                         n\sim N,m\sim M,k\sim K,\ell\sim L,
                         r\sim R,s\sim S,q\sim Q \\
                        1\leqslant N\leqslant M\leqslant K \leqslant L\leqslant y \\
                        1\leqslant R\leqslant S\leqslant Q\leqslant y
                            }}  g  \\
        &  &  \times  \min\left(T^{11/4},\frac{T^{9/4}}{|\sqrt{n}+\sqrt{m}+\sqrt{k}+\sqrt{\ell}
                       -\sqrt{r}-\sqrt{s}-\sqrt{q}| }\right).
\end{eqnarray*}
If $L\geqslant100Q,$ then $|\sqrt{n}+\sqrt{m}+\sqrt{k}+\sqrt{\ell}-\sqrt{r}-\sqrt{s}-\sqrt{q}|\gg L^{1/2},$ so the trivial estimate yields
\begin{equation*}
 G(N,M,K,L,R,S,Q)\ll\frac{T^{9/4+\varepsilon}NMKLRSQ}{(NMKLRSQ)^{3/4}L^{1/2}}\ll T^{9/4+\varepsilon}y^{5/4}\ll
  T^{41/16+\varepsilon}.
\end{equation*}
If $Q\geqslant100L,$ we can get the same estimate. So later we always suppose that $L\asymp Q.$ Let
$\eta=\sqrt{n}+\sqrt{m}+\sqrt{k}+\sqrt{\ell}-\sqrt{r}-\sqrt{s}-\sqrt{q}.$ Write
\begin{equation}\label{proposition-7}
   G(N,M,K,L,R,S,Q)=G_1+G_2+G_3,
\end{equation}
where
\begin{eqnarray*}
  & & G_1:=T^{11/4}\sum_{0<|\eta|\leqslant T^{-1/2}}g, \\
  & & G_2:=T^{9/4}\sum_{T^{-1/2}<|\eta|\leqslant 1}g|\eta|^{-1}, \\
  & & G_3:=T^{9/4}\sum_{|\eta|> 1}g|\eta|^{-1}.
\end{eqnarray*}
We estimate $G_1$ first. From $|\eta|\leqslant T^{-1/2},$ we get $Q\gg T^{1/63}$ via Lemma \ref{lemma-1}. By Lemma
\ref{lemma-8}, we get

\begin{eqnarray}\label{proposition-8}
   G_1 & \ll & \frac{T^{11/4+\varepsilon}}{(NMKLRSQ)^{3/4}}\mathscr{A}_1(N,M,K,L,R,S,Q;T^{-1/2})
                     \nonumber  \\
    & \ll & \frac{T^{11/4+\varepsilon}}{(NMKLRSQ)^{3/4}}\left(T^{-1/2}Q^{1/2}NMKLRS+NMKRSL^{1/2}\right)
                         \nonumber  \\
   & \ll & T^{9/4+\varepsilon}y^{5/4}+T^{11/4+\varepsilon}(NMKRS)^{1/4}L^{-1}
                         \nonumber  \\
   & \ll & T^{41/16+\varepsilon}+T^{11/4+\varepsilon}(NMKRS)^{1/4}L^{-1}.
\end{eqnarray}
By Lemma \ref{lemma-10}, we get
\begin{eqnarray}\label{proposition-9}
 G_1 & \ll & \frac{T^{11/4+\varepsilon}}{(NMKLRSQ)^{3/4}}\mathscr{A}_{-,-}(N,M,K,L,R,S,Q;T^{-1/2}) \nonumber\\
 & \ll & \frac{T^{11/4+\varepsilon}}{(NMKLRSQ)^{3/4}} (T^{-1/14}N^{13/14}+N^{5/7})(T^{-1/14}M^{13/14}+M^{5/7})
                   \nonumber \\
 &  & \times  (T^{-1/14}K^{13/14}+K^{5/7})
              (T^{-1/14}R^{13/14}+R^{5/7})(T^{-1/14}S^{13/14}+S^{5/7})
              (T^{-1/7}Q^{13/7}+Q^{10/7})
                   \nonumber \\
  & \ll & T^{11/4+\varepsilon} Q^{-1/14}(NMKRS)^{-1/28}(T^{-1/7}Q^{3/7}+1)(T^{-1/14}N^{3/14}+1) \nonumber  \\
  &  &    \times (T^{-1/14}M^{3/14}+1) (T^{-1/14}K^{3/14}+1) (T^{-1/14}R^{3/14}+1) (T^{-1/14}S^{3/14}+1)
                             \nonumber \\
 & \ll &  T^{11/4+\varepsilon}Q^{-1/14}(NMKRS)^{-1/28}
\end{eqnarray}
by noting that $T^{-1/7}D^{3/7}\ll1$ for $D=Q,N,M,K,R,S.$
From (\ref{proposition-8}) and (\ref{proposition-9}), we get
\begin{eqnarray}\label{proposition-10}
  G_1 & \ll & T^{41/16+\varepsilon}+T^{11/4+\varepsilon}\cdot\min\left(\frac{(NMKRS)^{1/4}}{L},
               \frac{1}{Q^{1/14}(NMKRS)^{1/28}} \right)
                    \nonumber \\
  & \ll & T^{41/16+\varepsilon}+T^{11/4+\varepsilon}\left(\frac{(NMKRS)^{1/4}}{L}\right)^{1/8}
          \left(\frac{1}{Q^{1/14}(NMKRS)^{1/28}} \right)^{7/8}
                     \nonumber \\
  & \ll & T^{41/16+\varepsilon}+T^{11/4+\varepsilon}Q^{-3/16}
                     \nonumber \\
  & \ll & T^{11/4-1/336+\varepsilon},
\end{eqnarray}
if we notice the fact that $Q\gg T^{1/63}.$

  Now we estimate $G_2.$ We also suppose $K\leqslant R$ and the other cases are the same. By a splitting argument,
we get the estimate
\begin{equation*}
   G_2\ll\frac{T^{9/4+\varepsilon}}{(NMKLRSQ)^{3/4}\delta}\times
          \sum_{\substack{\delta<|\eta|\leqslant2\delta\\|\eta|\neq0}}1
\end{equation*}
for some $T^{-1/2}\leqslant\delta\leqslant1.$ By Lemma \ref{lemma-8}, we get
\begin{eqnarray}\label{proposition-11}
   G_2 & \ll & \frac{T^{9/4+\varepsilon}}{(NMKLRSQ)^{3/4}\delta}\mathscr{A}_1(N,M,K,L,R,S,Q;2\delta)
                     \nonumber  \\
    & \ll & \frac{T^{9/4+\varepsilon}}{(NMKLRSQ)^{3/4}\delta}\left(\delta Q^{1/2}NMKLRS+NMKRSL^{1/2}\right)
                      \nonumber  \\
   & \ll & T^{9/4+\varepsilon}y^{5/4}+T^{9/4+\varepsilon}(NMKRS)^{1/4}(L\delta)^{-1}
                         \nonumber  \\
   & \ll & T^{41/16+\varepsilon}+T^{9/4+\varepsilon}(NMKRS)^{1/4}(L\delta)^{-1}.
\end{eqnarray}
On the other hand, by Lemma \ref{lemma-10}, we have
\begin{eqnarray*}
 G_2 & \ll & \frac{T^{9/4+\varepsilon}}{(NMKLRSQ)^{3/4}\delta}\mathscr{A}_{-,-}(N,M,K,L,R,S,Q;2\delta)
                 \nonumber \\
 & \ll & \frac{T^{9/4+\varepsilon}}{(NMKLRSQ)^{3/4}\delta} (\delta^{1/7}N^{13/14}+N^{5/7})(\delta^{1/7}M^{13/14}+M^{5/7})
                 \nonumber   \\
 &  & \times  (\delta^{1/7}K^{13/14}+K^{5/7})
              (\delta^{1/7}R^{13/14}+R^{5/7})(\delta^{1/7}S^{13/14}+S^{5/7})(\delta^{2/7}Q^{13/7}+Q^{10/7})
                 \nonumber   \\
  & \ll & T^{9/4+\varepsilon}\delta^{-1}
          (\delta^{1/7}N^{5/28}+N^{-1/28})(\delta^{1/7}M^{5/28}+M^{-1/28})(\delta^{1/7}K^{5/28}+K^{-1/28})
                              \nonumber \\
  &  & \times (\delta^{1/7}R^{5/28}+R^{-1/28})(\delta^{1/7}S^{5/28}+S^{-1/28})(\delta^{2/7}Q^{5/14}+Q^{-1/14})
                             \nonumber \\
  & \ll & T^{9/4+\varepsilon}\delta^{-1} Q^{-1/14}(NMKRS)^{-1/28}(\delta^{2/7}Q^{3/7}+1)(\delta^{1/7}N^{3/14}+1)
                               \nonumber  \\
  &  &    \times (\delta^{1/7}M^{3/14}+1) (\delta^{1/7}K^{3/14}+1) (\delta^{1/7}R^{3/14}+1) (\delta^{1/7}S^{3/14}+1)
                               \nonumber \\
  & \ll & T^{9/4+\varepsilon}\delta^{-1} Q^{-1/14}(NMKRS)^{-1/28}(\delta^{2/7}Q^{3/7}+1) (\delta^{5/7}(NMKRS)^{3/14}
                               \nonumber \\
  &  &    +\delta^{4/7}(MKRS)^{3/14}+\delta^{3/7}(KRS)^{3/14}+\delta^{2/7}(RS)^{3/14} +\delta^{1/7}(S)^{3/14}+1)
                               \nonumber \\
 & \ll & T^{9/4+\varepsilon}\delta^{-2/7} Q^{-1/14}(NMKRS)^{5/28}(\delta^{2/7}Q^{3/7}+1)
                               \nonumber \\
 &  &  +T^{9/4+\varepsilon}\delta^{-1}Q^{-1/14}(NMKRS)^{-1/28}(\delta^{2/7}Q^{3/7}+1)
                               \nonumber \\
 &  &  \times (\delta^{4/7}(MKRS)^{3/14}+\delta^{3/7}(KRS)^{3/14}+\delta^{2/7}(RS)^{3/14} +\delta^{1/7}S^{3/14}+1)
                               \nonumber \\
\end{eqnarray*}
\begin{eqnarray}\label{proposition-12}
 & \ll & T^{9/4+\varepsilon}y^{5/4}+T^{9/4+\varepsilon}T^{1/7}y^{23/28}+T^{9/4+\varepsilon}\delta^{-1}Q^{-1/14}
         (NMKRS)^{-1/28}(\delta^{2/7}Q^{3/7}+1)
                                \nonumber \\
   &   &    \times (\delta^{4/7}Q^{6/7}+\delta^{3/7}Q^{9/14}+\delta^{2/7}Q^{3/7}+\delta^{1/7}Q^{3/14}+1)
                                \nonumber \\
 & \ll &  T^{291/112+\varepsilon}+T^{9/4+\varepsilon}\delta^{-1}Q^{-1/14}(NMKRS)^{-1/28}(\delta^{6/7}Q^{9/7}+1).
\end{eqnarray}
From (\ref{proposition-11}) and (\ref{proposition-12}), we get
\begin{eqnarray}\label{proposition-13}
  G_2 & \ll & T^{291/112+\varepsilon}+T^{9/4+\varepsilon}\delta^{-1}\cdot\min\left(\frac{(NMKRS)^{1/4}}{L},
               \frac{\delta^{6/7}Q^{9/7}+1}{Q^{1/14}(NMKRS)^{1/28}} \right)
                         \nonumber \\
  & \ll & T^{291/112+\varepsilon}+T^{9/4+\varepsilon}\delta^{-1} \left(\frac{(NMKRS)^{1/4}}{L}\right)^{1/8}
                   \left(\frac{\delta^{6/7}Q^{9/7}+1}{Q^{1/14}(NMKRS)^{1/28}}\right)^{7/8}
                          \nonumber \\
 & \ll & T^{291/112+\varepsilon}+T^{9/4+\varepsilon}\delta^{-1}Q^{-3/16}(\delta^{3/4}Q^{9/8}+1).
\end{eqnarray}
If $\delta\gg Q^{-3/2},$ then $\delta^{3/4}Q^{9/8}\gg1$ and (\ref{proposition-13}) implies (recall $\delta\gg T^{-1/2}$)
\begin{eqnarray}\label{proposition-14}
  G_2 & \ll & T^{291/112+\varepsilon}+T^{9/4+\varepsilon}\delta^{-1/4}Q^{15/16}    \nonumber \\
 & \ll & T^{291/112+\varepsilon}+T^{9/4+1/8+\varepsilon}y^{15/16}  \nonumber \\
 & \ll & T^{167/64+\varepsilon}.
\end{eqnarray}
If $\delta\ll Q^{-3/2},$ then $\delta^{3/4}Q^{9/8}\ll1$ and (\ref{proposition-13}) implies
\begin{eqnarray}\label{proposition-15}
   G_2\ll T^{291/112+\varepsilon}+T^{9/4+\varepsilon}\delta^{-1}Q^{-3/16}.
\end{eqnarray}
By Lemma \ref{lemma-1}, we have $2\delta\geqslant|\eta|\gg Q^{-63/2},$ which implies $\delta^{-1}\ll Q^{63/2}.$
Therefore, we get
\begin{equation*}
  \delta^{-1}\ll\min\left(Q^{63/2},T^{1/2}\right).
\end{equation*}
From the above estimate and (\ref{proposition-15}), we get
\begin{eqnarray}\label{proposition-16}
   G_2 & \ll & T^{291/112+\varepsilon}+\min\left(T^{9/4+\varepsilon}Q^{501/16},T^{11/4+\varepsilon}Q^{-3/16}\right)
                        \nonumber \\
   & \ll & T^{291/112+\varepsilon}+\left(T^{9/4+\varepsilon}Q^{501/16}\right)^{1/168}
           \left(T^{11/4+\varepsilon}Q^{-3/16}\right)^{167/168}
                        \nonumber \\
   & \ll & T^{11/4-1/336+\varepsilon}.
\end{eqnarray}
For $G_3,$ by a splitting argument and Lemma \ref{lemma-8} (notice $|\eta|\gg1$),we get
\begin{eqnarray}\label{proposition-17}
  G_3 & \ll & \frac{T^{9/4+\varepsilon}}{(NMKLRSQ)^{3/4}\delta}\times
           \sum_{\substack{\delta<|\eta|\leqslant2\delta\\ |\delta|\gg1}}1
                           \nonumber \\
  & \ll & \frac{T^{9/4+\varepsilon}}{(NMKLRSQ)^{3/4}} Q^{1/2}NMKLRS
                             \nonumber \\
  & \ll & T^{9/4+\varepsilon}y^{5/4}
                             \nonumber \\
  & \ll &  T^{41/16+\varepsilon}.
\end{eqnarray}
Combining (\ref{proposition-6}), (\ref{proposition-7}), (\ref{proposition-10}), (\ref{proposition-16}) and (\ref{proposition-17}), we get
\begin{equation}\label{proposition-18}
  \int_T^{2T}S_2(x)\mathrm{d}x\ll T^{11/4-1/336+\varepsilon}.
\end{equation}
In the same way, we can prove that
\begin{equation}\label{proposition-19}
  \int_T^{2T}S_4(x)\mathrm{d}x\ll T^{11/4-1/336+\varepsilon}
\end{equation}
and
\begin{equation}\label{proposition-20}
  \int_T^{2T}S_6(x)\mathrm{d}x\ll T^{11/4-1/336+\varepsilon},
\end{equation}
if we use Lemma \ref{lemma-9} instead of Lemma \ref{lemma-8}.

Now Proposition \ref{mingti-1}, i.e. equation (\ref{proposition-1}), follows from (\ref{proposition-2})-(\ref{proposition-5}) and (\ref{proposition-18})-(\ref{proposition-20}).

\begin{proposition}\label{mingti-2}
   For $T^{\varepsilon}\ll y\ll T^{1/4},$ we have
 \begin{equation}\label{proposition-2-1}
   \int_T^{2T}|R_2|^7\mathrm{d}x\ll T^{11/4-1/336+\varepsilon}
 \end{equation}
 and
  \begin{equation}\label{proposition-2-2}
   \int_T^{2T}|R_1|^6|R_2|\mathrm{d}x\ll T^{11/4-1/336+\varepsilon}.
 \end{equation}
\end{proposition}

 \textbf{Proof.} First, for $T\leqslant x\leqslant2T,$ we have
\begin{eqnarray*}
   R_2 & = & \frac{1}{\sqrt{2}\pi}\sum_{y<n\leqslant T}\frac{d(n)}{n^{3/4}}x^{1/4}\cos(4\pi\sqrt{nx}-\pi/4)
              +O(x^{1/2+\varepsilon}T^{-1/2}) \\
  & \ll & x^{1/4}\left|\sum_{y<n\leqslant T}\frac{d(n)}{n^{3/4}}\cos(4\pi\sqrt{nx}-\pi/4)\right|+T^{\varepsilon} \\
 & \ll & x^{1/4}\left|\sum_{y<n\leqslant T}\frac{d(n)}{n^{3/4}}e\left(2\sqrt{nx}\right)\right|+T^{\varepsilon}.
\end{eqnarray*}
Therefore, one has
\begin{eqnarray}\label{proposition-2-3}
  \int_T^{2T}|R_2|^2\mathrm{d}x
 & \ll & T^{1+\varepsilon}+\int_T^{2T}x^{1/2}\left|\sum_{y<n\leqslant
                                    T}\frac{d(n)}{n^{3/4}}e\left(2\sqrt{nx}\right)\right|^2\mathrm{d}x
  \nonumber \\
  & \ll & T^{1+\varepsilon}+T^{1/2}\int_T^{2T}\sum_{y<n\leqslant T}\sum_{y<m\leqslant T}\frac{d(n)d(m)}{(nm)^{3/4}}
          e\left(2(\sqrt{n}-\sqrt{m})\sqrt{x}\right)\mathrm{d}x
                                               \nonumber \\
  & = & T^{1+\varepsilon}+T^{1/2}\sum_{y<n\leqslant T}\sum_{y<m\leqslant T}\frac{d(n)d(m)}{(nm)^{3/4}}
         \int_T^{2T}e\left(2(\sqrt{n}-\sqrt{m})\sqrt{x}\right)\mathrm{d}x
                                               \nonumber \\
  & \ll & T^{1+\varepsilon}+T^{3/2}\sum_{y<n\leqslant T}\frac{d^2(n)}{n^{3/2}}
          +T\mathop{\sum_{y<n\leqslant T}\sum_{y<m\leqslant T}}_{n\neq
                     m}\frac{d(n)d(m)}{(nm)^{3/4}|\sqrt{n}-\sqrt{m}|}
                                                \nonumber \\
  & \ll & T^{1+\varepsilon}+T^{3/2}\sum_{n>y}\frac{d^2(n)}{n^{3/2}}
          +T\sum_{m< n\leqslant T}\frac{d(n)d(m)}{(nm)^{3/4}|\sqrt{n}-\sqrt{m}|}
                                                 \nonumber \\
  & \ll & T^{1+\varepsilon}+T^{3/2}y^{-1/2}\log^3 y
         +T\sum_{m< n\leqslant T}\frac{d(n)d(m)}{(nm)^{3/4}(\sqrt{n}-\sqrt{m})}
                                                 \nonumber \\
  & \ll & T^{3/2}y^{-1/2}\log^3 T +T\sum_{m< n\leqslant T}\frac{d(n)d(m)}{(nm)^{3/4}(\sqrt{n}-\sqrt{m})}.
\end{eqnarray}
For the last sum in (\ref{proposition-2-3}), we have
\begin{equation}\label{proposition-2-4}
  \sum_{m< n\leqslant T}\frac{d(n)d(m)}{(nm)^{3/4}(\sqrt{n}-\sqrt{m})}:=\Sigma_1+\Sigma_2,
\end{equation}
where
\begin{eqnarray*}
  &   & \Sigma_1:=\sum_{\substack{m< n\leqslant T\\ \sqrt{n}-\sqrt{m}>\frac{1}{10}(nm)^{1/4}}}
           \frac{d(n)d(m)}{(nm)^{3/4}(\sqrt{n}-\sqrt{m})},   \nonumber \\
  &   & \Sigma_2:=\sum_{\substack{m< n\leqslant T\\ 0<\sqrt{n}-\sqrt{m}\leqslant\frac{1}{10}(nm)^{1/4}}}
           \frac{d(n)d(m)}{(nm)^{3/4}(\sqrt{n}-\sqrt{m})}.
\end{eqnarray*}
For $\Sigma_1,$ one has
\begin{equation}\label{proposition-2-5}
  \Sigma_1\ll\sum_{n\leqslant T}\sum_{m<n}\frac{d(n)d(m)}{nm}\ll\sum_{n\leqslant T}\frac{d(n)}{n}
              \sum_{m<n}\frac{d(m)}{m}\ll\log^4 T.
\end{equation}
For $\Sigma_2,$ we can write
\begin{equation}\label{proposition-2-6}
  \Sigma_2\ll\sum_{\substack{m< n\leqslant T\\ 0<\sqrt{n}-\sqrt{m}\leqslant\frac{1}{10}(nm)^{1/4}}}
           \frac{d^2(n)+d^2(m)}{(nm)^{3/4}(\sqrt{n}-\sqrt{m})}:=\Sigma_{21}+\Sigma_{22},
\end{equation}
where
\begin{eqnarray*}
  &   & \Sigma_{21}:=\sum_{\substack{m< n\leqslant T\\ 0<\sqrt{n}-\sqrt{m}\leqslant\frac{1}{10}(nm)^{1/4}}}
           \frac{d^2(n)}{(nm)^{3/4}(\sqrt{n}-\sqrt{m})},  \nonumber \\
  &   & \Sigma_{22}:=\sum_{\substack{m< n\leqslant T\\ 0<\sqrt{n}-\sqrt{m}\leqslant\frac{1}{10}(nm)^{1/4}}}
           \frac{d^2(m)}{(nm)^{3/4}(\sqrt{n}-\sqrt{m})}.
\end{eqnarray*}
By the mean value theorem and taking $n-m=r,$ we get
\begin{equation}\label{proposition-2-7}
 \Sigma_{21}\ll\sum_{n\leqslant T}\sum_{m<n}\frac{d^2(n)}{n(n-m)}\ll\sum_{n\leqslant T}\frac{d^2(n)}{n}
                \sum_{1\leqslant r\ll T}\frac{1}{r}\ll \log^5T,
\end{equation}
\begin{equation}\label{proposition-2-8}
  \Sigma_{22}\ll\sum_{n\leqslant T}\sum_{m<n}\frac{d^2(m)}{m(n-m)}\ll\sum_{m<T}\frac{d^2(m)}{m}
                \sum_{1\leqslant r\ll T}\frac{1}{r}\ll \log^5T.
\end{equation}
By the estimate (\ref{proposition-2-3})-(\ref{proposition-2-8}), we get
\begin{equation}\label{proposition-2-9}
    \int_T^{2T}|R_2|^2\mathrm{d}x\ll T^{3/2}y^{-1/2}\log^3 T.
\end{equation}

  Taking $y=T^{1/4},$ we get $R_1\ll T^{1/4}y^{1/4}\ll T^{131/416},$ and therefore $R_2\ll T^{131/416+\varepsilon}.$ Using Ivi\'{c}'s large-value technique directly to $R_2$ without modifications, one can get that the estimate
\begin{equation}\label{proposition-2-10}
   \int_T^{2T}|R_2|^{A_0}\mathrm{d}x\ll T^{1+\frac{A_0}{4}+\varepsilon}
\end{equation}
holds with $A_0=\frac{184}{19}, T^\varepsilon\ll y\ll T^{1/4}.$ Then, for $2<A<A_0,$ by (\ref{proposition-2-9}),
 (\ref{proposition-2-10}) and H\"{o}lder's inequality, we have
\begin{eqnarray*}
  \int_{T}^{2T}|R_2|^A\mathrm{d}x
       & = & \int_T^{2T} |R_2|^{\frac{2(A_0-A)}{A_0-2}+\frac{A_0(A-2)}{A_0-2}}\mathrm{d}x   \\
      & \ll & \left(\int_T^{2T}|R_2|^2 \mathrm{d}x \right)^{\frac{A_0-A}{A_0-2}}
              \left(\int_T^{2T}|R_2|^{A_0} \mathrm{d}x\right)^{\frac{A-2}{A_0-2}} \\
     & \ll & T^{1+\frac{A}{4}+\varepsilon}y^{-\frac{A_0-A}{2(A_0-2)}}.
\end{eqnarray*}
Especially, we take $A=7, y=T^{1/4}$ and get
\begin{eqnarray}\label{proposition-2-11}
   \int_T^{2T}|R_2|^7\mathrm{d}x & \ll & T^{11/4+\varepsilon}y^{-51/292}  \nonumber \\
   & \ll & T^{11/4-51/1168+\varepsilon} \nonumber \\
   & \ll & T^{11/4-1/336+\varepsilon}.
\end{eqnarray}

  By H\"{o}lder's inequality, we have
\begin{eqnarray}\label{proposition-2-12}
  \int_T^{2T}|R_1|^6|R_2|\mathrm{d}x
   & \ll & \left(\int_T^{2T}(|R_1|^6)^{7/6}\mathrm{d}x\right)^{6/7}
           \left(\int_T^{2T}|R_2|^7\mathrm{d}x\right)^{1/7}  \nonumber \\
   & \ll & \left(\int_T^{2T}|R_1|^7\mathrm{d}x\right)^{6/7}
           \left(\int_T^{2T}|R_2|^7\mathrm{d}x\right)^{1/7}  \nonumber \\
   & \ll & (T^{11/4})^{6/7}\cdot(T^{11/4-51/1168+\varepsilon})^{1/7}   \nonumber \\
   & \ll & T^{11/4-1/336+\varepsilon}.
\end{eqnarray}
Now Proposition \ref{mingti-2} follows from (\ref{proposition-2-11}) and (\ref{proposition-2-12}).

 From Proposition \ref{mingti-1} and Proposition \ref{mingti-2}, we get the seventh-power moment of $\Delta(x).$

\bigskip
\bigskip

\textbf{Acknowledgement}

   The author would like to express the most sincere gratitude to Professor Wenguang Zhai for his valuable
advice and constant encouragement.


\begin{thebibliography}{99}
  \bibitem {Corradi} K. Corr\'{a}di and I. K\'{a}tai, \textit{A comment on K. S. Gangadharan's paper entitled
   ``Two classical lattice point problems''}, Magyar Tud. Akad. Mat. Fiz. Oszt. Kozl. 17 (1967) 89-97.

  \bibitem {Hafner} J. L. Hafner, \textit{New omega theorems for two classical lattice point problems}, Invent.
                    Math. 63 (1981) 181-186.

  \bibitem {Hardy} G. H. Hardy, \textit{On Dirichlet's divisor problem}, Proc. London Math. Soc. (2) 15 (1916) 1-25.

                    \textit{Exponential sums with a difference}, Proc. London Math. Soc. (3) 61 (1990) 227-250.

  \bibitem {Heath-Brown} D. R. Heath-Brown, \textit{The distribution and moments of the error term in the Dirichlet
                            divisor problem}, Acta Arith. 60 (1992) 389-415.

                     66 (1993), 279-301.

  \bibitem {Huxley-2} M. N. Huxley, \textit{Exponential sums and lattice points III}, Proc. London Math. Soc. (3)
                     87 (2003) 591-609.

  \bibitem {Ivic} A. Ivi\'{c}, \textit{Large values of the error term in the divisor problem}, Invent. Math.
                     71 (1983) 513-520.
  \bibitem {Ivic-book} A. Ivi\'{c}, \textit{Lectures on mean values of the Riemann zeta-function}, Tata Institute
                    of Fundamental Research Lectures on Mathematics and Physics, Vol. 82, Bombay, 1991.

  \bibitem {Ivic-Sargos} A. Ivi\'{c} and P. Sargos, \textit{On the higher power moments of the error term in
                         the divisor problem}, Illinois J. Math. 51 (2007) 353-377.


  \bibitem {Kong} K. L. Kong, \textit{Some mean value theorems for certain error terms in analytic number theory},
             Master degree Thesis, The University of Hong Kong, 2014.

  \bibitem {Kratzel} E. Kr\"{a}tzel, \textit{Lattice points}, Deutsch. Verlag Wiss., Berlin, 1988.

  \bibitem {Piatetski-Shapiro} I. I. Piatetski-Shapiro, \textit{On a variant of Waring-Goldbach¡¯s problem}, Mat.
                        Sb. 30 (72) (1) (1952), 105-120 (in Russian).


  \bibitem {Robert-Sargos} O. Robert and P. Sargos, \textit{Three-dimensional exponential sums with monomials}, J.
                            Reine Angew. Math. 591 (2006) 1-20.

  \bibitem {Segal} B. I. Segal, \textit{On a theorem analogous to Waring¡¯s theorem}, Dokl. Akad. Nauk SSSR
                       (N. S.) 2 (1933), 47-49 (in Russian).

  \bibitem {Tong} K. C. Tong, \textit{On divisor problem III}, Ada Math. Sinica 6 (1956) 515-541.

  \bibitem {Tsang} K. M. Tsang, \textit{Higher-power moments of $\Delta(x),E(t)$ and $P(x)$}, Proc. London Math.
                    Soc. (3) 65 (1992) 65-84.

  \bibitem {Voronoi} G. F. Vorono\"{i}, \textit{Sur une fonction transcendante et ses applications \`{a} la sommation de quelques s\'{e}ries}, Ann. \'{E}cole Normale (3) 21 (1904) 207-268, 459-534.

  \bibitem {Wang} J. Wang, \textit{On the sixth-power moment of $\Delta(x)$ (in Chinese)}, Master degree Thesis,
                    Shandong Normal University, 2010.

  \bibitem {Zhai-1} W. G. Zhai, \textit{On higher-power moments of $\Delta(x)$}, Acta Arith. 112 (2004) 367-395.
  \bibitem {Zhai-2} W. G. Zhai, \textit{On higher-power moments of $\Delta(x)$ (II)}, Acta Arith. 114 (2004) 35-54.
  \bibitem {Zhai-3} W. G. Zhai, \textit{On higher-power moments of $\Delta(x)$ (III)}, Acta Arith.
                     118 (2005) 263-281.

  \bibitem {Zhang-deyu} D. Y. Zhang and W. G. Zhai, \textit{On the fifth-power moment of $\Delta(x)$}, International
    Journal of Number Theory, Vol. 7, No. 1 (2011) 71-86.

\end{thebibliography}
\end{document}